\documentclass[12pt]{article}

\usepackage{array,float,epsf,latexsym,psfrag,amssymb,eurosym,
theorem,amsmath,verbatim,framed}

\usepackage{graphicx}
\usepackage{oldgerm}
\usepackage{amsfonts}
\usepackage{bbm,mathbbol}
\usepackage{mathrsfs} 

\usepackage[dvips]{color}

\def\widebar   #1{\overline{#1}}

\def\el #1{\,{\mathbbm #1}}
\def\elb #1{\,\overline{\mathbbm #1}}
\def\elz #1{\,{\mathbb #1}}

\def\tz #1{\mbox{\boldmath $#1$}}

\def\gz #1{\mbox{\boldmath $\mit #1$}}

\def\gzb #1{\overline{\mbox{\boldmath $\mit #1$}}}

\def\tvsf #1{\mbox{\boldmath $\sf #1$}}

\def\und  {\quad\mbox{and }\quad}
\def\with {\quad\mbox{with}\quad}

\def\d {\,\mbox{d}}

\def\be {\begin{equation}}
\def\ee {\end  {equation}}
\def\bea {\begin{eqnarray}}
\def\eea {\end  {eqnarray}}
\def\bean {\begin{eqnarray*}}
\def\eean {\end  {eqnarray*}}

\def\b {{}^\bullet}
\def\db {{}^{\bullet\bullet}}
\def\c {{}^\circ}
\def\dc {{}^{\circ\circ}}

\begin{document}
\begin{center}
{\Large {\bf Analytical Mechanics Allows Novel Vistas on\\[2mm] Mathematical Epidemic Dynamics Modelling\footnote{\today
}}}\\[3mm]
P Steinmann$^{1,2}$ \\[2mm]
$^1$ Institute of Applied Mechanics, Friedrich-Alexander Universität Erlangen-Nürnberg,\\ 91054 Erlangen, Germany\\[1mm]
$^2$ Glasgow Computational Engineering Centre, University of Glasgow,\\ G12 8QQ Glasgow, United Kingdom\\[6mm]

{\bf Abstract}
\end{center}

This contribution aims to shed light on mathematical epidemic dynamics modelling from the viewpoint of analytical mechanics.
To set the stage, it recasts the basic SIR model of mathematical epidemic dynamics in an analytical mechanics setting. Thereby, it considers two possible re-parameterizations of the basic SIR model. On the one hand, it is proposed to re-scale time, while on the other hand, to transform the coordinates, i.e.\ the independent variables. In both cases, Hamilton's equations in terms of a suited Hamiltonian as well as Hamilton's principle in terms of a suited Lagrangian are considered in minimal and extended phase and state space coordinates, respectively. The corresponding Legendre transformations relating the various options for the Hamiltonians and Lagrangians are detailed. Ultimately, this contribution expands on a multitude of novel vistas on mathematical epidemic dynamics modelling that emerge from the analytical mechanics viewpoint. As result, it is believed that interesting and relevant new research avenues open up when exploiting in depth the analogies between analytical mechanics and mathematical epidemic dynamics modelling.

\section{Motivation}
The global COVID-19 pandemic, with alleged outbreak by the end of 2019 in Wuhan, China \cite{Lu2020ST}, -- despite its devastating implications for health, economy, and society -- has in particular challenged modelling and simulation of mathematical epidemic dynamics. Political decision makers around the globe seek (or should seek) advise from scientist such as virologists, biologists, clinicians, economists, sociologist as well as modelers from different fields. Especially the latter are in the position to virtually simulate various scenarios based on well-founded assumptions in order to provide support and guidance for the difficult and momentous decisions of politicians, e.g.\ on lockdown measures and step-wise exit strategies thereof. Thus, the critical importance of modelling is clearly appreciated, and indeed, mathematical epidemic dynamics modelling is a well-established and mature field:\\

Traditional mathematical modelling of epidemic dynamics roots in the concept of \textit{Susceptible}, \textit{Infected}, and \textit{Recovered} (SIR) compartments as originally proposed in \cite{Kermack1927M}. Various modifications extend the classical SIR model to account for further compartments such as, e.g., \textit{Deceased} (SIRD model), \textit{Exposed} (SEIR model), and \textit{Quarantined} (SIQRD model), among many other sophisticated options, see \cite{Hethcote2000,Diekmann2012HB}.  Classical SIR-type compartment-based models are coupled ordinary-differential-equations (ODEs). Extending the ODE-based SIR-type modelling approach to integro-differential-equations (IDEs) allows to also consider the detailed course of the disease, e.g., delay due to incubation time and the infectious period, see \cite{Keimer2020P}. SIR-type models describe the temporal spread of infectious diseases for integral populations, thereby, however, neglecting the interconnectedness of spatially distributed geographic areas. Recently proposed multiple compartment models, e.g. the mcSIR model in \cite{Seroussi2019LPSY}, also allow consideration of the geographical spread of potentially multiple infectious virus strains within the population and its potentially multiple subgroups (e.g. age groups). Spatial network models, e.g.\ \cite{Balcan2009,Balcan2010,Pastor2015}, for example based on the SEIR model at each network node can qualitatively simulate the outbreak dynamics of infectious diseases and the impact of travel restrictions in geographical areas at the global (macro) scale such as China and the USA \cite{Peirlinck2020} or Europe \cite{Linka2020}. However, due to its stochastic nature and strong impact of socio-economic factors, modelling epidemic dynamics within geographical areas at the local (micro) scale requires the use of another modelling paradigm, i.e. rule-driven, agent-based models \cite{Rahmandad2008S}. Agent-based models allow for example studying the effect of various lockdown exit strategies on local geographical entities with only a comparatively small number of individuals (agents), see e.g. \cite{German2020}. Regardless of the modelling approach taken, quantitative predictions of epidemic dynamics remain challenging and critically require careful identification of model parameters from reliable databases, see e.g.\ \cite{Kergassner2020}.

\begin{center}
\textit{Given this mature state of affairs, why is a novel, alternative view on mathematical epidemic dynamics modelling justified at all?}\\
\end{center}

The answer is as follows: To date, modelers of complex mechanical systems and behavior have developed a versatile and extremely successful tool-set, including sophisticated analytical and in particular efficient and accurate computational methods. Examples are techniques to master severe non-linearities and couplings with non-mechanical fields, a multitude of multi-scale and homogenization modelling approaches as well as incorporating uncertainty quantification into modelling and simulation. Mathematical epidemic dynamics modelling can undoubtedly benefit largely from this accumulated expertise! In summary, it is therefore believed that first translating epidemic dynamics models into an analytical mechanics setting (related steps towards this aim may be found, e.g., in \cite{Militaru2013M,Ionescu2015MM,Seroussi2019LPSY}) and then, secondly, exploiting the analogy between the two approaches while utilizing the full tool-set of mechanical modelling, can provide novel vistas and unprecedented opportunities. The present contribution aims to sketch out a few of these perspectives and to encourage the mechanics community to offer its strong modelling expertise to possibly and hopefully help further improve epidemic dynamics modelling.

\section{Basic SIR Model}
Classical modelling of epidemics dynamics roots in the concept of \textit{Susceptible}, \textit{Infected} and \textit{Recovered} (SIR) compartments as originally proposed in \cite{Kermack1927M}.
The basic compartment-based SIR  model is the following set of two coupled ordinary differential equations (ODEs)
\be
\left[\begin{array}{c}
I\b\\
S\b
\end{array}\right]=
\left[\begin{array}{l}
\phantom{+}\beta SI-\gamma I\\
         - \beta SI
\end{array}\right].
\ee
Here, $I$ and $S$ denote the stock of individuals in the \textit{Infected} and the \textit{Susceptible} compartments, respectively, normalized by the size of the entire population. The notation for the derivative of a quantity with respect to ordinary time $t$ is
\be
\{\}\b:=\d_t\{\}.
\ee
The parameters $\beta$ and $\gamma$ are the infection and the recovery rate, respectively, with their ratio defining the \textit{basic reproduction number} $R_0:=\beta/\gamma$. Note finally that the stock of individuals in the \textit{Recovered} compartment follows from the constraint $S+I+R=1$, thus the evolution equation $R\b=-[S\b+I\b]$ is tacitly suppressed in our presentation.

\section{Time Re-Parameterized SIR Model}
In order to recast the basic SIR model into a format more amenable to the analytical mechanics setting, it is proposed, as a first option, to re-scale time as
\be
\tau:=\int_0^t SI\d\bar t\with\d\tau/\d t=SI.
\ee
The notation for the derivative of a quantity with respect to re-scaled time $\tau$ is
\be
\{\}\c:=\d_\tau\{\}.
\ee
Consequently, the derivatives with respect to re-scaled and ordinary time are related via
\be
\{\}\c=\{\}\b/[SI].
\ee
As a result, the basic SIR model is eventually re-parameterized in terms of re-scaled time as
\be\label{timeSIR}
\left[\begin{array}{c}
I\c\\
S\c
\end{array}\right]=\gz F\with\gz F:=
\left[\begin{array}{l}
\phantom{+}\beta-\gamma/S\\
         - \beta
\end{array}\right]=:
\left[\begin{array}{c}
 V\\
 F
\end{array}\right].
\ee
In the time re-parameterized SIR model the right-hand-side is abbreviated as the \textit{forcing term} $\gz F$, i.e.\ as the column matrix consisting of the time re-parameterized \textit{rate of infection} $V$ and \textit{force of infection} $F$ (a common terminology from mathematical epidemiology already establishing a semantic analogy to mechanics). Computing $I\dc=[\gamma/S^2]S\c$ from Eq.\ \ref{timeSIR}.1 and eliminating $S\c$ with the help of Eq.\ \ref{timeSIR}.2 renders $I\dc=-\beta\gamma/S^2$.
Resorting finally again to Eq.\ \ref{timeSIR}.1, i.e.\ expressing $S(I\c)=\gamma/[\beta-I\c]$, allows formulating the time re-parameterized SIR model as a single non-linear ODE
\be\label{timeSIRsingle}
I\dc=-R_0\,[\beta-I\c]^2
\ee
exclusively in the stock of individuals in the \textit{Infected} compartment and with the basic reproduction number and the infection rate as parameters.

\subsection{Hamiltonian in Minimal Phase Space Coordinates}
The minimal phase space coordinates collectively assembled in the column matrix $\gz Z\in\elz R^2$, i.e.\ the \textit{generalized coordinate} $Q$ defined as the stock of individuals in the \textit{Infected} compartment and the \textit{generalized momentum} $P$ defined as the stock of individuals in the \textit{Susceptible} compartment, span the two-dimensional phase space $\el P$, thus
\be
\el P:=\left\{\gz Z:=\left[\begin{array}{c}
Q\\
P
\end{array}\right]:=
\left[\begin{array}{c}
I\\
S
\end{array}\right]\right\}.
\ee
The Hamiltonian $H(\gz Z)$ in minimal phase space coordinates, which eventually results in the time re-parameterized SIR model from Eq.\ \ref{timeSIR}, is then identified as
\be\label{HamiltonianMinPhase}
H(\gz Z):=\beta[I+S]-\gamma\ln S.
\ee
Indeed, the corresponding Hamilton equations deliver a re-formulation of the result in Eq.\ \ref{timeSIR}, i.e.\
\be
\left[\begin{array}{c}
I\c\\
S\c
\end{array}\right]=
\left[\begin{array}{cc}
\phantom{-}0&1\\
-1&0
\end{array}\right]
\left[\begin{array}{l}
\beta\\
\beta-\gamma/S
\end{array}\right].
\ee
Symbolic notation clearly reveals the Hamiltonian structure of the time re-parameterized SIR model
\be\fbox{$\displaystyle
\gz Z\c=\gz J\cdot\partial_{\gz Z} H(\gz Z)$}
\ee
Thereby, the skew-symmetric $\gz J\in\elz R^2\times\elz R^2$ denotes the so-called \textit{symplectic matrix} in $\el P$
\be
\gz J:=
\left[\begin{array}{cc}
\phantom{-}0&1\\
-1&0
\end{array}\right]
\ee
with $\gz J^t=-\gz J$ and $\gz J^2=-\gz I$, whereby $\gz I\in\elz R^2\times\elz R^2$ is the common unit matrix in $\el P$. The Hamiltonian structure in terms of the skew-symmetric $\gz J$ clearly identifies the (autonomous) Hamiltonian $H(\gz Z)$ in minimal phase space coordinates as first integral, i.e.\ as a conserved quantity under the flow implied by the Hamilton equations, since
\be
H\c=\partial_{\gz Z} H\cdot\gz Z\c=\partial_{\gz Z} H\cdot\gz J\cdot\partial_{\gz Z} H=0.
\ee
The gradient of the Hamiltonian $H(\gz Z)$ with respect to the minimal phase space coordinates, abbreviated in the sequel as $\gz G(\gz Z)\in\elz R^2$, computes as
\be
\partial_{\gz Z} H(\gz Z)=
\left[\begin{array}{l}
\beta\\
\beta-\gamma/S
\end{array}\right]=:\gz G(\gz Z).
\ee
Taken together, the time re-parameterized SIR model obeys Hamiltonian structure and identifies the relation between the gradient $\gz G(\gz Z)\in\elz R^2$ of the Hamiltonian $H(\gz Z)$ in minimal phase space coordinates and the forcing term $\gz F(\gz Z)\in\elz R^2$ as
\be\label{ForceMinPhase}\fbox{$\displaystyle
\gz Z\c=\gz J\cdot\gz G(\gz Z)=\gz F(\gz Z)$}
\ee
Exchanging the time derivative to the one with respect to ordinary time $t$ destroys the clean Hamiltonian structure in terms of a constant symplectic matrix, i.e.\
\be
\gz Z\b=[SI]\,\gz J\cdot\gz G(\gz Z)=[SI]\,\gz F(\gz Z).
\ee
One may, of course, re-interpret this result as a Hamiltonian structure with non-constant, coordinate-dependent symplectic matrix $[SI]\,\gz J$ on a non-flat manifold.

\subsection{Lagrangian in Minimal State Space Coordinate}
The minimal state space coordinate, i.e.\ the \textit{generalized coordinate} $Q$ defined as the stock of individuals in the \textit{Infected} compartment, spans the one-dimensional state space $\el S$, thus
\be
\el S:=\left\{Q:=I\right\}.
\ee
Then a Legendre transformation of the Hamiltonian in minimal phase space coordinates defines the corresponding Lagrangian
\be
L(I,I\c):=\sup_S\{\,S\,I\c-H(\gz Z)\,\}.
\ee
The supremum condition identifies $I\c$ with the derivative $\partial_SH(\gz Z)$ of the Hamiltonian in minimal phase space coordinates from Eq.\ \ref{HamiltonianMinPhase} with respect to the generalized momentum, and renders
\be
I\c\doteq\beta-\gamma/S.
\ee
Re-solving the above supremum condition for $S$ in terms of $I\c$ delivers
\be
S(I\c)=\gamma/[\beta-I\c].
\ee
Then, with $S(I\c)\,I\c=\gamma\,I\c/[\beta-I\c]$ and $H\big(I, S(I\c)\big)=\beta\big[I+\gamma/[\beta-I\c]\big]-\gamma\ln(\gamma/[\beta-I\c])$, the Lagrangian in minimal state space coordinate follows eventually as
\be
L(I,I\c)=-\beta I-\gamma+\gamma\ln(\gamma/[\beta-I\c]).
\ee
Based on the Lagrangian in minimal state space coordinates, Hamilton's principle results in the stationarity condition
\be
\left[\frac{\partial L}{\partial I\c}\right]^\circ=\frac{\partial L}{\partial I}.
\ee
Thus, with $\partial_{I\c}L=\gamma/[\beta-I\c]\,\to\,[\partial_{I\c}L]\c=\gamma\,I\dc/[\beta-I\c]^2$ and $\partial_IL=-\beta$, the Euler-Lagrange equation corresponding to the Lagrangian in minimal state space coordinates reads
\be\fbox{$\displaystyle
I\dc+R_0\,[\beta-I\c]^2=0$}
\ee
Clearly, the Euler-Lagrange equation in minimal state space coordinates coincides with the single, non-linear ODE formulation of the time re-parameterized SIR model in Eq.\ \ref{timeSIRsingle}.

\subsection{Lagrangian in Extended State Space Coordinates}
Alternatively, extended state space coordinates collectively assembled in the column matrix $\gz Q\in\elz R^2$, i.e.\ the \textit{generalized coordinates} jointly defined as the stock of individuals in the \textit{Infected} and \textit{Susceptible} compartments, span the two-dimensional state space $\elb S$, thus
\be
\elb S:=\left\{\gz Q:=\left[\begin{array}{c} I\\ S\end{array}\right]\right\}.
\ee
Then the Lagrangian $\widebar L(\gz Q, \gz Q\c)$ in extended state space coordinates, which eventually results in the time re-parameterized SIR model from Eq.\ \ref{timeSIR}, is determined as
\be\label{LagrangianExtPhase}
\widebar L(\gz Q, \gz Q\c):=\frac{1}{2}\gz Q\cdot\gz J{}^t\cdot\gz Q\c-H(\gz Q).
\ee
Here, the Legendre transformation term $\gz Q\cdot\gz J^t\cdot\gz Q\c=\gz Q\c\cdot\gz J\cdot\gz Q$ expands as the skew-symmetric form
\be
\gz Q\cdot\gz J{}^t\cdot\gz Q\c=S\,I\c-I\,S\c,
\ee
whereas the Hamiltonian $H(\gz Q)$, corresponding to Eq. \ref{HamiltonianMinPhase}, is now parameterized in terms of the extended state space coordinates $\gz Q$ with
\be
H(\gz Q):=\beta[I+S]-\gamma\ln S.
\ee
As a result, the Lagrangian in extended state space coordinates then renders the stationarity conditions of the corresponding Hamilton's principle as
\be
\left[\frac{\partial\widebar L}{\partial\gz Q\c}\right]^\circ=\frac{\partial\widebar L}{\partial\gz Q}.
\ee
With $\partial_{\gz Q\c}\widebar L=\gz Q\cdot\gz J^t/2$ and $\partial_{\gz Q}\widebar L=\gz Q\c\cdot\gz J/2-\gz G(\gz Q)$, the Euler-Lagrange equations thus follow as
\be\label{ELextState}
\frac{1}{2}\gz Q\c\cdot\gz J^t=\frac{1}{2}\gz Q\c\cdot\gz J-\gz G(\gz Q).
\ee
Unfolding the compact symbolic notation, these expand concretely into
\be
\frac{1}{2}\left[\begin{array}{c}
\phantom{+}S\c\\
         - I\c
\end{array}\right]=
\frac{1}{2}
\left[\begin{array}{c}
         - S\c\\
\phantom{+}I\c
\end{array}\right]-
\left[\begin{array}{l}
\beta\\
\beta-\gamma/S
\end{array}\right].
\ee
The Euler-Lagrange equations in Eq.\ \ref{ELextState} are next re-formulated by recalling that due to the skew-symmetry of the symplectic matrix, $\gz Q\c\cdot\gz J=-\gz Q\c\cdot\gz J^t$ and $\gz Q\c\cdot\gz J^t=\gz J\cdot\gz Q\c$ hold, thus
\be\label{ELextStateRef}
\gz J\cdot\gz Q\c=-\gz G(\gz Q).
\ee
Finally, with $\gz J^2=-\gz I$, the re-formulated Euler-Lagrange equations in Eq.\ \ref{ELextStateRef} are modified into
\be\fbox{$\displaystyle
\gz Q\c=\gz J\cdot\gz G(\gz Q)=\gz F(\gz Q)$}
\ee
Obviously, this format recovers the relation between the gradient $\gz G(\gz Q)\in\elz R^2$ of the Hamiltonian $H(\gz Q)$ (in extended state space coordinates) and the forcing term $\gz F(\gz Q)\in\elz R^2$ already established previously in Eq. \ref{ForceMinPhase}.
Likewise, exchanging the time derivative to the one with respect to ordinary time $t$
\be
\gz Q\b=[SI]\,\gz J\cdot\gz G(\gz Q)=[SI]\,\gz F(\gz Q)
\ee
results again in a formulation with non-constant, coordinate dependent symplectic matrix $[SI]\,\gz J$ on a non-flat manifold.

\subsection{Hamiltonian in Extended Phase Space Coordinates}
For academic curiosity it is also interesting to consider extended phase space coordinates $\gz Q$, i.e.\ the \textit{generalized coordinates} jointly defined as the stock of individuals in the \textit{Infected} and \textit{Susceptible} compartments, and $\gz P$, i.e.\ heretofore undefined \textit{generalized momenta}, collectively assembled in the column matrix $\gzb Z\in\elz R^4$, to
span the four-dimensional state space $\elb P$, thus
\be{}
\elb P:=\left\{\gzb Z:=\left[\begin{array}{c}\gz Q\\\gz P\end{array}\right]\with
\gz Q:=\left[\begin{array}{c}I\\ S\end{array}\right]\und
\gz P:=\left[\begin{array}{c}\mathit{\Upsilon}\\\mathit{\Sigma}\end{array}\right]\right\}.
\ee
A Legendre transformation of the Lagrangian in extended state space coordinates from Eq. \ref{LagrangianExtPhase} defines the associated Hamiltonian
\be
\widebar H(\gz Q,\gz P)=\sup_{\gz Q\c}\{\, \gz P\cdot\gz Q\c-\widebar L(\gz Q, \gz Q\c)\,\}.
\ee
The corresponding supremum condition identifies $\gz P$ with the derivative $\partial_{\gz Q\c}\widebar L(\gz Q, \gz Q\c)$ of the Lagrangian $\widebar L(\gz Q, \gz Q\c)$ in extended state space coordinates from Eq.\ \ref{LagrangianExtPhase} with respect to the time re-parameterized rate of the generalized coordinates $\gz Q\c$, and renders
\be
\frac{\partial\widebar L}{\partial\gz Q\c}=\frac{1}{2}\gz Q\cdot\gz J{}^t=:\gz P.
\ee
As a result, $\gz P$ does not depend on $\gz Q\c$, thus identifying the Lagrangian $\widebar L(\gz Q, \gz Q\c)$ in Eq.\ \ref{LagrangianExtPhase} as \textit{degenerate} in the sense of the Dirac's generalized Hamiltonian dynamics \cite{Dirac1950}. Consequently, with $\gz Q\cdot\gz J^t=\gz J\cdot\gz Q$ and $\gz J^2=-\gz I$, an additional constraint for the extended phase space coordinates emerges
\be
\gz C(\gz Q,\gz P)=\gz Q+2\gz J\cdot\gz P\doteq\tz 0.
\ee
The Hamiltonian $\widebar H(\gz Q,\gz P)$ in extended phase space coordinates thus follows from Legendre transformation by incorporating the constraint via the Lagrange multiplier $\gz\Lambda$, i.e.\
\be
\widebar H(\gz Q,\gz P)=\frac{1}{2}\gz Q\cdot\gz J{}^t\cdot\gz Q\c-\widebar L(\gz Q, \gz Q\c)+\gz\Lambda\cdot\gz C(\gz Q,\gz P).
\ee
Taking into account the explicit form of the Lagrangian $\widebar L(\gz Q, \gz Q\c)$ in extended state space coordinates from Eq.\ \ref{LagrangianExtPhase} then renders the explicit representation of the Hamiltonian $\widebar H(\gz Q,\gz P)$ in extended phase space coordinates
\be\label{HamiltonianExtPhase}
\widebar H(\gz Q,\gz P)=
H(\gz Q)+\gz\Lambda\cdot\gz C(\gz Q,\gz P).
\ee
Invoking $\partial_{\gz Q}\gz C(\gz Q,\gz P)=\gz I$ and $\partial_{\gz P}\gz C(\gz Q,\gz P)=2\gz J$, Hamilton's equations based on the Hamiltonian $\widebar H(\gz Q,\gz P)$ in extended phase space coordinates from Eq. \ref{HamiltonianExtPhase} result in
\be
\left[\begin{array}{c}
\gz Q\c\\
\gz P\c
\end{array}\right]=
\left[\begin{array}{cc}
\phantom{-}0&\gz I\\
-\gz I&0
\end{array}\right]
\left[\begin{array}{l}
\partial_{\gz Q}\widebar H\\
\partial_{\gz P}\widebar H
\end{array}\right]=
\left[\begin{array}{cc}
\phantom{-}0&\gz I\\
-\gz I&0
\end{array}\right]
\left[\begin{array}{c}
\gz G(\gz Q)+\gz\Lambda\\
2\gz\Lambda\cdot\gz J
\end{array}\right].
\ee
The unknown Lagrange multiplier $\gz\Lambda$ is determined from the consistency condition for the constraint
\be
\gz C\c(\gz Q\c,\gz P\c)=\gz Q\c+2\gz J\cdot\gz P\c\doteq\tz 0,
\ee
which, upon introducing Hamilton's equations $\gz Q\c=2\gz\Lambda\cdot\gz J$ and $\gz P\c=-[\gz G(\gz Q)+\gz\Lambda]$, results in
\be
2\gz\Lambda\cdot\gz J-2\gz J\cdot[\gz G(\gz Q)+\gz\Lambda]\doteq\tz 0.
\ee
The solution of the consistency condition then renders the explicit representation for the Lagrange multiplier
\be
\gz\Lambda=-\frac{1}{2}\gz G(\gz Q).
\ee
Using again Hamilton's equation $\gz Q\c=2\gz\Lambda\cdot\gz J$ and exploiting the skew-symmetry of the symplectic matrix, i.e.\ $\gz\Lambda\cdot\gz J=-\gz J\cdot\gz\Lambda$, recovers once more  the already previously established relation between the gradient $\gz G(\gz Q)\in\elz R^2$ of the Hamiltonian $H(\gz Q)$ (in extended state space coordinates) and the forcing term $\gz F(\gz Q)\in\elz R^2$
\be\fbox{$\displaystyle
\gz Q\c=\gz J\cdot\gz G(\gz Q)=\gz F(\gz Q)$}
\ee
Finally, using Hamilton's equation $\gz P\c=-[\gz G(\gz Q)+\gz\Lambda]$ identifies eventually the time re-parameterized rate of the heretofore unknown generalized momenta as
\be\fbox{$\displaystyle
\gz P\c=-\frac{1}{2}\gz G(\gz Q)=\frac{1}{2}\gz J\cdot\gz F(\gz Q)$}
\ee
Thereby, the last equality follows from $\gz J^2=-\gz I$ and the previous relation $\gz J\cdot\gz G(\gz Q)=\gz F(\gz Q)$.

\section{Coordinate Re-Parameterized SIR Model}
Alternatively, and again aiming to recast the basic SIR model into an analytical mechanics format, it is proposed, as a second option, to logarithmically transform the coordinates (or rather the independent variables) as
\be
I\mapsto i:=\ln I\und
S\mapsto s:=\ln S.
\ee
As a result, the basic SIR model, re-parameterized in terms of logarithmically transformed coordinates, however still in terms of derivatives with respect to ordinary time $t$, now reads
\be\label{coorSIR}
\left[\begin{array}{c}
i\b\\
s\b
\end{array}\right]=\gz f
\with
\gz f:=
\left[\begin{array}{l}
\phantom{+}\beta   S( s)-\gamma\\
         - \beta\, I( i)
\end{array}\right]=:
\left[\begin{array}{c}
v\\
f
\end{array}\right].
\ee
Note that in this re-parametrization, the stocks of individuals in the \textit{Infected} and \textit{Susceptible} compartments are considered as \textit{dependent functions} of the logarithmically transformed coordinates
\be
I( i):=\exp i\und
S( s):=\exp s.
\ee
In the coordinate re-parameterized SIR model the right-hand-side is abbreviated as the \textit{forcing term} $\gz f$, i.e.\ as the column matrix consisting of the coordinate re-parameterized \textit{rate of infection} $v$ and \textit{force of infection} $f$. Computing $i\db=\beta\,S(s)\,s\b$ from Eq.\ \ref{coorSIR}.1 and eliminating $s\b$ with the help of Eq.\ \ref{coorSIR}.2 renders $i\db=-\beta^2S(s)\,I(i)$.
Resorting finally again to Eq.\ \ref{coorSIR}.1, i.e.\ expressing $S(s)=[i\b+\gamma]/\beta$, allows formulating the coordinate re-parameterized SIR model as a single non-linear ODE
\be\label{coorSIRsingle}
i\db=-\beta\,I(i)\,[i\b+\gamma]
\ee
exclusively in the logarithmic stock of individuals in the \textit{Infected} compartment.

\subsection{Hamiltonian in Minimal Phase Space Coordinates}
The minimal (logarithmic) phase space coordinates collectively assembled in the column matrix $\gz z\in\elz R^2$, i.e.\ the \textit{generalized coordinate} $q$ defined as the \textit{logarithmic} stock of individuals in the \textit{Infected} compartment and the \textit{generalized momentum} $p$ defined as the \textit{logarithmic} stock of individuals in the \textit{Susceptible} compartment, span the two-dimensional phase space $\el p$, thus
\be
\el p:=\left\{\gz z:=\left[\begin{array}{c}
q\\
p
\end{array}\right]:=
\left[\begin{array}{c}
 i\\
 s
\end{array}\right]\right\}.
\ee
Next, the Hamiltonian $h(\gz z)$ in minimal (logarithmic) phase space coordinates, which eventually results in the coordinate re-parameterized SIR model from Eq.\ \ref{coorSIR}, is identified as
\be\label{HamiltonianMinPhaseCoor}
h(\gz z):=\beta[I(i)+S(s)]-\gamma s.
\ee
As a result, the corresponding Hamilton equations deliver a re-formulation of Eq.\ \ref{coorSIR}, i.e.\
\be
\left[\begin{array}{c}
i\b\\
s\b
\end{array}\right]=
\left[\begin{array}{cc}
\phantom{-}0&1\\
-1&0
\end{array}\right]
\left[\begin{array}{l}
\beta\, I( i)\\
\beta   S( s)-\gamma
\end{array}\right].
\ee
Symbolic notation showcases clearly the Hamiltonian structure of the (logarithmic) coordinate re-parameterized SIR model
\be\fbox{$\displaystyle
\gz z\b=\gz j\cdot\partial_{\gz z} h(\gz z)$}
\ee
Thereby, the skew-symmetric $\gz j\in\elz R^2\times\elz R^2$ denotes the appropriate symplectic matrix in $\el p$
\be
\gz j:=
\left[\begin{array}{cc}
\phantom{-}0&1\\
-1&0
\end{array}\right]
\ee
with $\gz j^t=-\gz j$ and $\gz j^2=-\gz i$, whereby $\gz i\in\elz R^2\times\elz R^2$ is the common unit matrix in $\el p$. The Hamiltonian structure in terms of the skew-symmetric $\gz j$ clearly identifies the (autonomous) Hamiltonian $h(\gz z)$ in minimal (logarithmic) phase space coordinates as first integral
\be
h\b=\partial_{\gz z} h\cdot\gz z\b=\partial_{\gz z} h\cdot\gz j\cdot\partial_{\gz z} h=0.
\ee
The gradient of the Hamiltonian $h(\gz z)$ with respect to the minimal (logarithmic) phase space coordinates, abbreviated in the sequel as $\gz g(\gz z)\in\elz R^2$, computes as
\be
\partial_{\gz z} h(\gz z)=
\left[\begin{array}{l}
\beta\, I( i)\\
\beta   S( s)-\gamma
\end{array}\right]=:\gz g(\gz z).
\ee
Summarizing, the coordinate re-parameterized SIR model obeys Hamiltonian structure and identifies the relation between the gradient $\gz g(\gz z)\in\elz R^2$ of the Hamiltonian $h(\gz z)$ in minimal (logarithmic) phase space coordinates and the forcing term $\gz f(\gz z)\in\elz R^2$ as
\be\label{ForceMinPhaseCoor}\fbox{$\displaystyle
\gz z\b=\gz j\cdot\gz g(\gz z)=\gz f(\gz z)$}
\ee
It is emphasized that despite the logarithmic nature of the re-parameterized coordinates, here, as a benefit, the Hamiltonian structure involves the constant, coordinate-independent symplectic matrix $\gz j$ of a flat manifold as well as derivatives with respect to ordinary time $t$.

\subsection{Lagrangian in Minimal State Space Coordinate}
The minimal (logarithmic) state space coordinate, i.e.\ the \textit{generalized coordinate} $q$ defined as the \textit{logarithmic} stock of individuals in the \textit{Infected} compartment, span the one-dimensional state space $\el s$, thus
\be
\el s:=\left\{q:=i\right\}.
\ee
Legendre transformation of the Hamiltonian in minimal (logarithmic) phase space coordinates defines the corresponding Lagrangian
\be
l(i,i\b)=\sup_s\{\,s\,i\b-h(\gz z)\,\}.
\ee
Then the supremum condition identifies $i\b$ with the derivative $\partial_sh(\gz z)$ of the Hamiltonian in minimal (logarithmic) phase space coordinates from Eq.\ \ref{HamiltonianMinPhaseCoor} with respect to the generalized momentum, and renders
\be
i\b=\beta S(s)-\gamma.
\ee
Re-solving the above supremum condition for $s$ in terms of $i\b$ delivers
\be
s(i\b)=\ln([i\b+\gamma]/\beta).
\ee
Then, with $s(i\b)\,i\b=\ln([i\b+\gamma]/\beta)\,i\b$ and $h\big(i, s(i\b)\big)=\beta\big[I(i)+[i\b+\gamma]/\beta\big]-\gamma\,\ln([i\b+\gamma]/\beta)$, the Lagrangian in minimal (logarithmic) state space coordinate follows eventually as
\be
l(i,i\b)=-\beta I(i)-[i\b+\gamma][1-\ln([i\b+\gamma]/\beta)].
\ee
Based on the Lagrangian in minimal (logarithmic) state space coordinates, Hamilton's principle results in the stationarity condition
\be
\left[\frac{\partial l}{\partial i\b}\right]^\bullet=\frac{\partial l}{\partial i}.
\ee
Thus, with $\partial_{i\b}l=-[1-\ln([i\b+\gamma]/\beta)]+\beta\,\to\,[\partial_{i\b}l]\b=i\db/[i\b+\gamma]$ and $\partial_il=-\beta\,I(i)$, the Euler-Lagrange equation corresponding to the Lagrangian in minimal (logarithmic) state space coordinates reads
\be\fbox{$\displaystyle
i\db+\beta\,I(i)\,[i\b+\gamma]=0$}
\ee
Obviously, the Euler-Lagrange equation in minimal (logarithmic) state space coordinates coincides with the single, non-linear ODE formulation of the coordinate re-parameterized SIR model in Eq.~\ref{coorSIRsingle}.

\subsection{Lagrangian in Extended State Space Coordinates}
Alternatively, extended (logarithmic) state space coordinates collectively assembled in the column matrix $\gz q\in\elz R^2$, i.e.\ the \textit{generalized coordinates} jointly defined as the \textit{logarithmic} stock of individuals in the \textit{Infected} and \textit{Susceptible} compartments, span the two-dimensional state space $\elb s$, thus
\be
\elb s:=\left\{\gz q:=\left[\begin{array}{c}i\\ s\end{array}\right]\right\}.
\ee
The Lagrangian $\widebar l(\gz q, \gz q\b)$ in extended (logarithmic) state space coordinates, which eventually results in the coordinate re-parameterized SIR model from Eq.\ \ref{coorSIR}, reads
\be\label{LagrangianExtPhaseCoor}
\widebar l(\gz q, \gz q\b):=\frac{1}{2}\gz q\cdot\gz j^t\cdot\gz q\b-h(\gz q).
\ee
Here, the Legendre transformation term $\gz q\cdot\gz j^t\cdot\gz q\b=\gz q\b\cdot\gz j\cdot\gz q$ expands as the skew-symmetric form
\be
\gz q\cdot\gz j^t\cdot\gz q\b=
s\,i\b-i\,s\b.
\ee
The Hamiltonian $h(\gz q)$, corresponding to Eq. \ref{HamiltonianMinPhaseCoor}, is now parameterized in terms of the extended (logarithmic) state space coordinates $\gz q$ as
\be
h(\gz q):=\beta[I(i)+S(s)]-\gamma s.
\ee
In extended (logarithmic) state space coordinates the Lagrangian then renders the stationarity conditions of the corresponding Hamilton's principle
\be
\left[\frac{\partial\widebar l}{\partial\gz q\b}\right]^\bullet=\frac{\partial\widebar l}{\partial\gz q}.
\ee
With $\partial_{\gz q\b}\widebar l=\gz q\cdot\gz j^t/2$ and $\partial_{\gz q}\widebar l=\gz q\b\cdot\gz j/2-\gz g(\gz q)$, the Euler-Lagrange equations thus follow as
\be\label{ELextStateCoor}
\frac{1}{2}\gz q\b\cdot\gz j^t=\frac{1}{2}\gz q\b\cdot\gz j-\gz g(\gz q)
\ee
and concretely unfold into
\be
\frac{1}{2}\left[\begin{array}{c}
\phantom{+}s\b\\
         - i\b
\end{array}\right]=
\frac{1}{2}
\left[\begin{array}{c}
         - s\b\\
\phantom{+}i\b
\end{array}\right]-
\left[\begin{array}{l}
\beta I( i)\\
\beta S( s)-\gamma
\end{array}\right].
\ee
Subsequently, the Euler-Lagrange equations in Eq.\ \ref{ELextStateCoor} are re-formulated by noting that $\gz q\b\cdot\gz j=-\gz q\b\cdot\gz j^t$ and $\gz q\b\cdot\gz j^t=\gz j\cdot\gz q\b$, thus
\be\label{ELextStateRefCoor}
\gz j\cdot\gz q\b=-\gz g(\gz q).
\ee
Finally, with $\gz j^2=-\gz i$, the Euler-Lagrange equations from Eq.\ \ref{ELextStateRefCoor} re-formulate as
\be\fbox{$\displaystyle
\gz q\b=\gz j\cdot\gz g(\gz q)=\gz f(\gz q)$}
\ee
This format recovers the relation between the gradient $\gz g(\gz q)\in\elz R^2$ of the Hamiltonian $h(\gz q)$ (in extended (logarithmic) state space coordinates) and the forcing term $\gz f(\gz q)\in\elz R^2$ already established previously in Eq. \ref{ForceMinPhaseCoor}.

\subsection{Hamiltonian in Extended Phase Space Coordinates}
For completeness we also consider extended (logarithmic) phase space coordinates $\gz q$, i.e.\ the \textit{generalized coordinates} jointly defined as the \textit{logarithmic} stock of individuals in the \textit{Infected} and \textit{Susceptible} compartments, and $\gz p$, i.e.\ heretofore undefined \textit{generalized momenta}, collectively assembled in the column matrix $\gzb z\in\elz R^4$, to
span the four-dimensional state space $\elb p$, thus
\be{}
\elb p:=\left\{\gzb z:=\left[\begin{array}{c}\gz q\\\gz p\end{array}\right]\with
\gz q:=\left[\begin{array}{c}i\\ s\end{array}\right]\und
\gz p:=\left[\begin{array}{c}\upsilon\\\sigma\end{array}\right]\right\}.
\ee
A Legendre transformation of the Lagrangian in extended (logarithmic) state space coordinates from Eq. \ref{LagrangianExtPhaseCoor} defines the associated Hamiltonian
\be
\widebar h(\gz q,\gz p)=\sup_{\gz q\b}\{\, \gz p\cdot\gz q\b-\widebar l(\gz q, \gz q\b)\,\}.
\ee
The Legendre transformation identifies $\gz p$ with the derivative $\partial_{\gz q\b}\widebar l(\gz q, \gz q\b)$ of the Lagrangian $\widebar l(\gz q, \gz q\b)$ in extended (logarithmic) state space coordinates from Eq.\ \ref{LagrangianExtPhaseCoor} with respect to the rate of the generalized coordinates $\gz q\b$, and renders
\be
\frac{\partial\widebar l}{\partial\gz q\b}=\frac{1}{2}\gz q\cdot\gz j^t=:\gz p.
\ee
Thus $\gz p$ does not depend on $\gz q\b$, consequently identifying the Lagrangian $\widebar l(\gz q, \gz q\b)$ in Eq.\ \ref{LagrangianExtPhaseCoor} as \textit{degenerate}, \cite{Dirac1950}. Consequently, with $\gz q\cdot\gz j^t=\gz j\cdot\gz q$ and $\gz j^2=-\gz i$, an additional constraint for the extended (logarithmic) phase space coordinates results
\be
\gz c(\gz q,\gz p)=\gz q+2\gz j\cdot\gz p\doteq\tz 0.
\ee
The Hamiltonian $\widebar h(\gz q,\gz p)$ in extended (logarithmic) phase space coordinates thus follows from Legendre transformation by incorporating the constraint via the Lagrange multiplier $\gz\lambda$, i.e.\
\be
\widebar h(\gz q,\gz p)=\frac{1}{2}\gz q\cdot\gz j^t\cdot\gz q\b-\widebar l(\gz q, \gz q\b)+\gz\lambda\cdot\gz c(\gz q,\gz p).
\ee
Taking into account the explicit form of the Lagrangian $\widebar l(\gz q, \gz q\b)$ in extended (logarithmic) state space coordinates from Eq.\ \ref{LagrangianExtPhaseCoor} results in the explicit representation of the Hamiltonian $\widebar h(\gz q,\gz p)$ in extended (logarithmic) phase space coordinates
\be\label{HamiltonianExtPhaseCoor}
\widebar h(\gz q,\gz p)=h(\gz q)+\gz\lambda\cdot\gz c(\gz q,\gz p).
\ee
Invoking $\partial_{\gz q}\gz c(\gz q,\gz p)=\gz i$ and $\partial_{\gz p}\gz c(\gz q,\gz p)=2\gz j$, Hamilton's equations based on the Hamiltonian $\widebar h(\gz q,\gz p)$ in extended (logarithmic) phase space coordinates from Eq. \ref{HamiltonianExtPhaseCoor} result in
\be
\left[\begin{array}{c}
\gz q\b\\
\gz p\b
\end{array}\right]=
\left[\begin{array}{cc}
\phantom{-}0&\gz i\\
-\gz i&0
\end{array}\right]
\left[\begin{array}{l}
\partial_{\gz q}\widebar h\\
\partial_{\gz p}\widebar h
\end{array}\right]=
\left[\begin{array}{cc}
\phantom{-}0&\gz i\\
-\gz i&0
\end{array}\right]
\left[\begin{array}{c}
\gz g(\gz q)+\gz\lambda\\
2\gz\lambda\cdot\gz j
\end{array}\right].
\ee
The unknown Lagrange multiplier $\gz\lambda$ is determined from the consistency condition for the constraint
\be
\gz c\b(\gz q\b,\gz p\b)=\gz q\b+2\gz j\cdot\gz p\b=\tz 0,
\ee
which, with $\gz q\b=2\gz\lambda\cdot\gz j$ and $\gz p\b=-[\gz g(\gz q)+\gz\lambda]$, results in
\be
2\gz\lambda\cdot\gz j-2\gz j\cdot[\gz g(\gz q)+\gz\lambda]\doteq\tz 0.
\ee
The solution of the consistency condition then renders the explicit representation for the Lagrange multiplier
\be
\gz\lambda=-\frac{1}{2}\gz g(\gz q).
\ee
Using again Hamilton's equation $\gz q\b=2\gz\lambda\cdot\gz j$ and exploiting $\gz\lambda\cdot\gz j=-\gz j\cdot\gz\lambda$, recovers once more  the already previously established relation between the gradient $\gz g(\gz q)\in\elz R^2$ of the Hamiltonian $h(\gz q)$ (in extended (logarithmic) state space coordinates) and the forcing term $\gz f(\gz q)\in\elz R^2$
\be\fbox{$\displaystyle
\gz q\b=\gz j\cdot\gz g(\gz q)=\gz f(\gz q)$}
\ee
Moreover, using Hamilton's equation $\gz p\b=-[\gz g(\gz q)+\gz\lambda]$ identifies finally the rate of the heretofore unknown generalized momenta as
\be\fbox{$\displaystyle
\gz p\b=-\frac{1}{2}\gz g(\gz q)=\frac{1}{2}\gz j\cdot\gz f(\gz q)$}
\ee
The last equality follows from $\gz j^2=-\gz i$ and the previous relation $\gz j\cdot\gz g(\gz q)=\gz f(\gz q)$.

\section{Novel Vistas from Analytical Mechanics}

In order to systematically explore lessons that can be learned from an analytical mechanics viewpoint on mathematical epidemic dynamics modelling, the coordinate re-parameterized version of the SIR model, based on the Hamiltonian in phase space, is taken as the point of departure:

\be\fbox{$\displaystyle
\gz z\b=\gz j\cdot\gz g(\gz z)$}
\ee

Building on this compact representation, several - analytical-mechanics-inspired - novel vistas on mathematical epidemic dynamics that promise fruitful research avenues for its modelling are identified  in the sequel.

\begin{description}
\item[Vista 1:] Allowing for non-autonomous, i.e.\ time-dependent generalized Hamiltonians results in the generalized representation
\be
\gz z\b=\gz j\cdot\gz g(\gz z, t).
\ee
Possible options justifying a non-autonomous Hamiltonian are for example:
\begin{description}
\item[Option a)] Various lockdown measures (cancellation of large events, school closing, contact limitations, etc) as well as their reversal (exit strategies) at discrete points in time are modelled by time-dependent parameters such as for example the infection and the recovery rates
    $$\beta=\beta(t)\und\gamma=\gamma(t),$$
    thus making the Hamiltonian time-dependent.
\item[Option b)] Various modifications extend the classical SIR model to account for  further compartments such as, e.g., \textit{Deceased} (SIRD model), \textit{Exposed} (SEIR model), \textit{Quarantined} (SIQRD model), among many other, more sophisticated options, see \cite{Hethcote2000,Diekmann2012HB}. SIR+ models of these types are then captured by appropriately extending the phase space variables
    $$\gz z:=[I,S,\cdots]\in\elz R^{2+\cdots}$$
    contribution to the Hamiltonian and its gradient.
\item[Option c)] Classical SIR-type compartment based models are coupled ordinary differential equations (ODEs). Extending the ODE-based SIR-type modelling approach to integro differential equations allows also considering, e.g., delay due to incubation time and infectious period, see \cite{Keimer2020P}. For $p$ representing relevant parameters, e.g.\ continuously distributed risk groups and/or past time, from space $P$, the right-hand-side of these read as
    $$\gz g(\gz z, t)=\int_P\gz j\cdot\gz\gamma(\gz z, p, t)\d p$$
    with $\gz\gamma(\gz z, p, t)$ the appropriate $p$-density of $\gz g(\gz z, t)$.
\end{description}
\item[Vista 2:] Allowing for an infinite dimensional phase space with its coordinates $\gz z=\gz z(\gz x, t)\in\elz R^{2+\cdots}$ defined as fields in four-dimensional space-time, results in the generalized representation
\be
\gz z\b=\gz j\cdot\gz g\big\{\gz z(\gz x),t\big\}.
\ee
Here, the right-hand-side is a functional of the phase space coordinates rather than a function.
Possible options for an infinite dimensional phase space are for example:
\begin{description}
\item[Option a)] Gradient-type models, whereby the Hamiltonian depends on the phase space coordinates $\gz z=\gz z(\gz x, t)$ and their higher spatial gradients
    $$\gz g\big\{\gz z(\gz x),t\big\}=\delta_{\gz z} h(\gz z(\gz x),\nabla_{\gz x}\gz z(\gz x),\cdots, t).$$
    Consequently, the right-hand-side follows from the variational derivative of the Hamiltonian, rather than from its gradient. Partial differential equations of reaction-convection-diffusion-type describing the spatio-temporal spread of infectious diseases are thus a modelling option \cite{Yamazaki2017W}.
\item[Option b)] Integral-type models, whereby, similar to Peridynamics formulations \cite{Javili2020MS}, the right-hand-side follows from a spatial integration
    $$\gz g\big\{\gz z(\gz x),t\big\}=\displaystyle\int_X\gz\gamma\big(\gz z(\gz x),\gz x, t\big)\d\gz x$$ over a cut-off domain $X$ (horizon)  that covers spatial interaction.
\end{description}
\item[Vista 3:] Allowing for an finite dimensional phase space with its coordinates $\tvsf z\in\elz R^{[2+\cdots]n}$ defined as column matrices, results in the generalized representation
\be
\tvsf z\b=\tvsf j\cdot\tvsf g(\tvsf z, t).
\ee
Possible options for an finite dimensional phase space are for example:
\begin{description}
\item[Option a)] Partition of the entire population into sub-populations, thereby separately considering different age/gender/risk groups, see e.g.\ \cite{Pastor2015}
    $$\tvsf z=[\gz z_{{\rm pop}_1},\cdots,\gz z_{{\rm pop}_{max}}].$$
\item[Option b)] Partition into various geographical locations in general network models accounting for the spatio-temporal spread of infectious diseases, see e.g.\ \cite{Seroussi2019LPSY}
    $$\tvsf z=[\gz z_{{\rm loc}_1},\cdots,\gz z_{{\rm loc}_{max}}].$$
\item[Option c)] Partition into multiple virus strains providing for generic infectious diseases, see e.g.\ \cite{Levy2018IY}
    $$\tvsf z=[\gz z_{{\rm vir}_1},\cdots,\gz z_{{\rm vir}_{max}}].$$
\end{description}
Capturing the interactions among the various partitions in a network is then reflected by the off-diagonal terms in the Hessian $\tvsf H:=\partial^2_{\tvsf z\tvsf z}h$ of the Hamiltonian.
\item[Vista 4:] For a pandemic such as COVID-19, spatial (geographical) resolution, i.e.\ resolution of a network is required at multiple scales: at the global (macro) scale, i.e.\ for the entire globe; at the medium (meso) scale, i.e.\ for individual countries; and at the local (micro) scale, i.e.\ for individual cities/communities. A fully detailed spatial resolution at the local (micro) scale for the entire globe is computationally prohibitive, moreover, most often an overkill degree of detail is also not needed and/or not possible due to the lack of data. However, the spatial resolution shall be adaptive to the quantity of interest, e.g.\ to study the dynamics of infectious disease spread in a particular city/community, only the integral results of more remote locations on the globe matter. These can be captured by a reduced resolution of the network in that geographical remote locations. This asks for a true multi-scale approach that adaptively zooms in only where needed.
    Possible options for multi-scaling are for example:
\begin{description}
\item[Option a)] Vertical coupling of scales relies on the assumption that the two scales considered are sufficiently separated, see e.g.\ \cite{Saeb2016SJ}. Then the 'force' term on the right-hand-side can be up-scaled from a sub-scale model by averaging in the sense of computational homogenisation
    $$\gzb z\b=<\tvsf j\cdot\tvsf g(\tvsf z, t;\gzb z)>.$$
    Here, the sup-scale (indicated by an over-bar) at the left-hand-side behaves like a SIR-type model whereas the sub-scale model at the right-hand-side lives either on a finite dimensional phase space or is represented by a rule-driven, so-called agent-based model. Agent-based models are an alternative modelling paradigm considering only a comparatively small number of individuals (agents). They are capable to capture the stochastic nature and strong impact of socio-economic factors present at small scales, see, e.g., \cite{German2020,Rahmandad2008S}. The sub-scale model is driven by the sup-scale phase space coordinates, whereby a proper scale-transition condition defines suited boundary/initial conditions at the sub-scale.
\item[Option b)] Horizontal coupling of scales, analogously to the Quasi-Continuum method \cite{Miller2002T}, requires adaptive resolution of the network spacing, here indicated by the sup-script $h$
    $$\gz z\b_h=\tvsf j\cdot\tvsf g(\tvsf z_h, t).$$
    Adaptivity requires suited network densification indicators that may follow from a proper error analysis, a topic that is still largely under-investigated for epidemic dynamics models.
\end{description}
\item[Vista 5:] The availability and reliability of recorded data, e.g.\ regarding the cumulative or daily infection cases, during an epidemic is typically characterized by a large degree of uncertainty, e.g.\ regarding the infection rates, the degree of immunity and/or their dark figures. Uncertainty quantification is based on simulations with uncertain data
\be
\tvsf z\b=\tvsf j\cdot\tvsf g\big(\tvsf z(\omega), t\big).
\ee
Thereby, uncertain data is here parameterized in terms of elementary events $\omega$ from which one may repeatedly draw samples to investigate uncertainty propagation throughout our model, see e.g.\ \cite{Pivovarov2018OWS,Pivovarov2019WS}.
Possible options for the description of uncertainties are for example:
\begin{description}
\item[Option a)] Aleatoric uncertainties require the use of random variables with probability density function (pdf) as a measure of likelihood  (e.g.\ Gaussian pdf in terms of the mean value and standard deviation). Aleatoric uncertainties are stochastic by nature and may not be neglected when the standard deviation is large.
\item[Option b)] Epistemic uncertainties may be captured by fuzzy variables with possibility density function as a measure of degree-of-membership (e.g.\ symmetric triangular membership function in terms of its modal value and support). Epistemic uncertainties reflect a lack of knowledge and, in the case of epidemic dynamics modelling, can be reduced by increasing testing for either infections and/or for anti-bodies.
\end{description}
\item[Vista 6:] The discrete trajectory in time of the phase space variables is algorithmically traced by an integrator of the generic format
\be
\tvsf z^{n+1}=\tvsf z^n+\Delta t\,\tvsf j\cdot\tvsf g\big(\tvsf z^{n+\alpha}(\omega), t^{n+\alpha}\big).
\ee
Here, sub-scripts $n+1$, $n$, and $n+\alpha$ refer to, respectively, the end point, the start point, and an intermediate point of/within a time step of length $\Delta t$. Possible options for time integrators that display different accuracy, stability, and robustness, in particular when integrating non-linear right-hand-sides are for example:
\begin{description}
\item[Option a)] Runge-Kutta integrators are of-the-shelf algorithms that come in a variety of different flavours (following from the corresponding Butcher tableau) like, e.g., single- and multi-stage integrators of varying algorithmic accuracy. However, they do not necessarily respect first integrals such as the conservation of the Hamiltonian for autonomous cases and may thus suffer from long-term deterioration of algorithmic stability and robustness.
\item[Option b)] Variational integrators are based on a discrete form of the action integral, whereby the integrand of the action integral is given by a discrete Lagrangian
    $$\int_{t^n}^{t^{n+1}}l(\tvsf z^{n+1},\tvsf z^n,\Delta t)\d t\to\mbox{stat}.$$
    The resulting discrete Hamiltonian principle then renders a variational integrator that follows from the discrete action integral being stationary. Variational integrators preserve symmetries (momentum maps) and structure (symplecticity) and are thus characterized by long-term algorithmic accuracy, stability, and robustness, see e.g.\ \cite{Lew2004MOW}.
\item[Option c)] Time-Finite-Element integrators follow from discretizing the Galerkin (weak) form of the Hamilton equations
    $$\int_T\delta\tvsf z\cdot[\tvsf z\b-\tvsf j\cdot\tvsf g]\d t=0\quad\forall\delta\tvsf z.$$
    Choosing appropriate Ansatz spaces for the test and trial functions, and suited quadrature rules for approximating the time integrals render integrators of arbitrary algorithmic accuracy that are also characterized by long-term algorithmic stability and robustness, see, e.g.\ \cite{Betsch2000S}.
\end{description}
\item[Vista 7:] The underlying equations governing epidemic dynamics are oftentimes unknown. However, they may be discovered from a data-driven approach \cite{Brunton2016} if sufficient many data is available, a scenario that is typically met for the spatio-temporal spread of infectious diseases. The key idea is then to connect the matrix arrangement of available discrete data points for the rate of the phase space coordinates by a matrix of Ansatz functions, e.g.\ monomials of the matrix arrangement of available discrete data points for the phase space coordinates, with the matrix arrangement of discrete Ansatz parameters.
\be
[\tvsf z\b_{{\rm dat}_1},\cdots,\tvsf z\b_{{\rm dat}_{max}}]^{\sf T}=\tvsf A(\tvsf z_{{\rm dat}_1}, \cdots,\tvsf z_{{\rm dat}_{max}})[\tvsf a_{{\rm par}_1}, \cdots,\tvsf a_{{\rm par}_{max}}]^{\sf T}.
\ee
Only few relevant entries in the matrix arrangement of the discrete Ansatz parameters are then determined from sparse regression, consequently the resulting models are denoted as parsimonious and compromise between accuracy and complexity, while avoiding overfitting.
\item[Vista 8:] Many more modelling approaches inspired by analytical mechanics are conceivable and it is left to the mechanics community to harness those to further improve mathematical epidemic dynamics modelling.
\end{description}

\section{Summary and Perspective}

First, this contribution explored options of how to recast the classical SIR model of mathematical epidemic dynamics modelling in the variational setting of analytical mechanics. In particular, it demonstrated that two conceptually entirely different re-parameterizations of the basic SIR model, i.e.\ either by re-scaling time or by transforming coordinates (independent variables), severely ease identification of corresponding Hamiltonians and Lagrangians for use within Hamilton's equations and Hamilton's principle. In each case, formulations in either minimal or extended phase and state space coordinates are possibly, providing in total eight different modelling options. Interestingly, in minimal phase space coordinates, the stock of individuals in the \textit{Infected} and the \textit{Susceptible} compartments represent the generalized coordinate \textit{and} the generalized momentum, respectively. In contrast, for extended phase space coordinates, they \textit{jointly} represent the generalized coordinates, whereas the associated generalized momenta are \textit{initially unknown} and only follow from exploiting a constraint on the extended phase space coordinates. However, regardless of the particular formulation chosen, from either Hamilton's equations or Hamilton's principle one eventually recovers the original set of coupled ODEs of the SIR model. As a recommendation, logarithmically transforming the coordinates appears more attractive, since derivatives with respect to ordinary time are retained for the evolution of the phase space coordinates. As an important perspective, recasting the classical SIR model in one of the eight different modelling options enables the analytical mechanician to employ the full mechanical modelling tool-set for a plethora of important extensions. The striking analogy between analytical mechanics and mathematical epidemic dynamics modelling opens up a multitude of fascinating and relevant new research avenues for the progression of the latter. It is thus believed that future exploitation of the Hamiltonian and/or Lagrangian structure of mathematical epidemic dynamics modelling leads to unprecedented insights and options for novel formulations.

\end{document}